\documentclass[12pt,twoside]{amsart}
\usepackage{amsfonts,amsmath}
\usepackage[all,knot,arc]{xy}

\def\endpf{\relax\ifmmode\expandafter\endproofmath\else
  \unskip\nobreak\hfil\penalty50\hskip.75em\hbox{}\nobreak\hfil\bull
  {\parfillskip=0pt \finalhyphendemerits=0 \bigbreak}\fi}
\def\bull{\vbox{\hrule\hbox{\vrule\kern3pt\vbox{\kern6pt}\kern3pt\vrule}\hrule}}

\newtheorem{defn}{Definition}[section]
\newtheorem{lemma}[defn]{Lemma}

\newtheorem{remark}[defn]{Remark}
\newtheorem{proposition}[defn]{Proposition}

\newtheorem{maintheorem}{Theorem}
\newtheorem{maincor}[maintheorem]{Corollary}
\newtheorem{maindef}[maintheorem]{Definition}

\newcommand{\cc}{{\mathbb C}}
\newcommand{\zz}{{\mathbb Z}}
\newcommand{\rr}{{\mathbb R}}

\newcommand{\qq}{{\mathbb Q}}
\newcommand{\pp}{{\mathbb P}}

\newcommand{\ozsvath}{Ozsv\'{a}th}
\newcommand{\szabo}{Szab\'{o}}

\newcommand{\spin}{\ifmmode{\rm Spin}\else{${\rm spin}$\ }\fi}
\newcommand{\spinc}{\ifmmode{{\rm Spin}^c}\else{${\rm spin}^c$\ }\fi}

\newcommand{\spincs}{\mathfrak s}

\newcommand{\ds}{\displaystyle}

\newenvironment{narrow}[2]{%
 \begin{list}{}{%
  \setlength{\topsep}{0pt}%
  \setlength{\leftmargin}{#1}%
  \setlength{\rightmargin}{#2}%
  \setlength{\listparindent}{\parindent}%
  \setlength{\itemindent}{\parindent}%
  \setlength{\parsep}{\parskip}%
 }%
\item[]}{\end{list}}

\newif\ifpic
\picfalse    
\pictrue   

\DeclareMathOperator\diag{Diag}

\DeclareMathOperator\Span{Span}

\begin{document}

\title{On slicing invariants of knots}
\author{Brendan Owens}
\date{\today}
\thanks{B. Owens was supported in part by NSF grant DMS-0604876.}

\begin{abstract}
The slicing number of a knot, $u_s(K)$, is the minimum number of crossing changes required
to convert $K$ to a slice knot.   This invariant is bounded above by the unknotting number and below
by the slice genus $g_s(K)$.  We show that for many knots, previous bounds on unknotting number 
obtained by \ozsvath\ and \szabo\ and by the author in fact give bounds on the slicing number.
Livingston defined another invariant $U_s(K)$ which takes into account signs of crossings changed
to get a slice knot, and which is bounded above by the slicing number and below by the slice genus.
We exhibit an infinite family of knots $K_n$ with slice genus $n$ and Livingston invariant greater than $n$.
Our bounds are based on restrictions (using Donaldson's diagonalisation
theorem or Heegaard
Floer homology) on the intersection forms of four-manifolds bounded by the
double branched cover of a knot.
\end{abstract}

\maketitle

\pagestyle{myheadings} \markboth{BRENDAN OWENS}
{ON SLICING INVARIANTS OF KNOTS}


\section{Introduction}
\label{sec:intro}
The unknotting number of a knot is the minimum number of crossing changes required to convert it to an unknot.  \ozsvath\ and \szabo\ used Heegaard Floer theory to provide a powerful obstruction to a knot having unknotting number one \cite{osu1}.  This obstruction was generalised in \cite{u} to higher unknotting numbers.  In this paper we show that similar techniques yield information about the number of crossing changes required to convert to a slice knot.

The slice genus $g_s(K)$ of a knot $K$ in the three-sphere is the minimum genus
of a connected oriented smoothly properly embedded surface in the four-ball with
boundary $K$.  A knot is called \emph{slice} if $g_s(K)=0$.
Given any diagram $D$ for a knot $K$, a new
knot may be obtained by changing one or more crossings of $D$. The
slicing number $u_s(K)$ is the minimum number of crossing changes
required to obtain a slice knot, where the minimum is taken over all
diagrams for $K$.  A ``movie" of a sequence of crossing changes represents
an immersed annulus in $S^3\times [0,1]$ with a singularity for each crossing change.
A neighbourhood of each singular point may be removed and replaced with an annulus;
if the last frame of the movie is a slice knot it may be capped off with a disk, yielding
a surface in $B^4$ with genus $u_s(K)$ and boundary $K$.

\begin{figure}[tbp] 
\begin{center}
\ifpic
\leavevmode
\begin{xy}
0;/r2pc/:
(0,0)*{
\begin{xy}
0;/r6pc/:
(0,0)*{}="1";
(1,0)*{}="2";
(1,1.2)*{}="3";
(0,1.2)*{}="4";
"1";"3" **\crv{}?(1)*\dir{>}; \POS?(.5)*{\hole}="x"; 
"2";"x" **\crv{}; 
"x";"4" **\crv{}?(1)*\dir{>}; 
(0.5,-0.2)*{{\rm Positive}}; 
\end{xy}};
(5,0)*{
\begin{xy}
0;/r6pc/:
(0,0)*{}="1";
(1,0)*{}="2";
(1,1.2)*{}="3";
(0,1.2)*{}="4";
"2";"4" **\crv{}?(1)*\dir{>}; \POS?(.5)*{\hole}="x"; 
"1";"x" **\crv{}; 
"x";"3" **\crv{}?(1)*\dir{>}; 
(0.5,-0.2)*{{\rm Negative}}; 
\end{xy}}
\end{xy}
\else \vskip 5cm \fi
\begin{narrow}{0.3in}{0.3in}
\caption{
\bf{Signed crossings in a knot diagram.}}
\label{fig:crossings}
\end{narrow}
\end{center}
\end{figure}
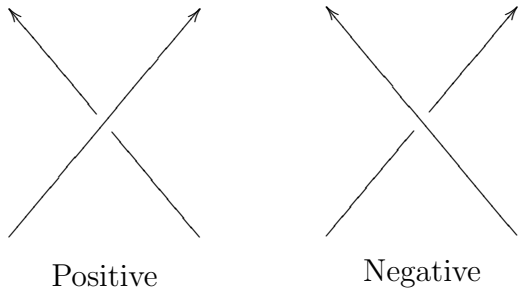

Recall that crossings in a knot diagram may be given a sign as in
Figure \ref{fig:crossings} (independent of the choice of
orientation of the knot).  Suppose that $K$ may be \emph{sliced} 
(converted to a slice knot) by changing $p$ positive and
$n$ negative crossings (in some diagram).  Form the immersed annulus in
$S^3\times [0,1]$ as before.  The sign of each self-intersection of this annulus
agrees with the sign of the corresponding crossing in the changed diagram.  Take two
self-intersections of opposite sign, and in each case remove a disk neighbourhood of
the singularity from just one of the intersecting sheets and connect the boundary components
by a tube.  This leads to a surface in $B^4$ with genus $\max(p,n)$.  Livingston defined
the following slicing invariant:
$$U_s(K)=\min(\max(p,n)),$$
where the minimum is taken over all diagrams for $K$ and over all sets of crossing changes
in a diagram which give a slice knot.

From the preceding discussion we see that
$$g_s(K)\le U_s(K)\le u_s(K).$$
Livingston showed in \cite{liv} that the two-bridge knot
$S(15,4)$, also known as $7_4$, has $u_s=2$ and $g_s=1$, thus giving a negative
answer to a question of Askitas \cite{ask}.  (Murakami and Yasuhara showed in \cite{my} 
that $8_{16}$
has $u_s=2$ and $g_s=1$.  Their proof is based on a four-manifold bounded by the double-branched cover of the knot.  We take a similar approach here.)  
Livingston also asked whether in fact $U_s$ is always equal to the slice genus, and suggested that $7_4$ may be a counterexample.

Let $\sigma(K)$ denote the signature of
a knot $K$. It is shown in \cite[Proposition 2.1]{cl} (also \cite[Theorem 5.1]{st}) 
that if $K'$ is obtained from
$K$ by changing a positive crossing, then
$$\sigma(K')\in\{\sigma(K),\sigma(K)+2\};$$
similarly if $K'$ is obtained from $K$ by changing a negative
crossing then
$$\sigma(K')\in\{\sigma(K),\sigma(K)-2\}.$$
Now suppose that $K$ may be sliced by changing $p$ positive and
$n$ negative crossings (in some diagram).  Since a slice knot has
zero signature, it follows that a bound for $n$ is given by
\begin{equation}
\label{eqn:nsig} n\ge\sigma(K)/2. \end{equation}
In this paper we give an obstruction to equality in
(\ref{eqn:nsig}).

Let $\Sigma(K)$ denote the double cover of $S^3$ branched along
$K$, and suppose that crossing changes in some diagram for $K$ result in a slice knot $J$.
It follows from ``Montesinos' trick'' (\cite{mont}, or see \cite{u}) that $\Sigma(K)$
is given by Dehn surgery on some framed link in $\Sigma(J)$
with half-integral framing coefficients.  

\begin{maindef}
\label{def:halfint}
An integer-valued symmetric bilinear form $Q$ on a free abelian group of rank $2r$ is 
said to be of \emph{half-integer surgery type}
if it admits a basis $\{x_1,\dots,x_r,y_1,\dots,y_r\}$ with
\begin{eqnarray*}
Q(x_i,x_j)&=&2\delta_{ij},\\
Q(x_i,y_j)&=&\delta_{ij}.
\end{eqnarray*}
\end{maindef}

\noindent{\bf Examples.} The positive-definite rank 2 unimodular form 
is of half-integer surgery type
since it may be represented by the matrix $\left(\begin{matrix}2&1\\1&1\end{matrix}\right)$.
The form represented by the matrix $\left(\begin{matrix}4&1\\1&4\end{matrix}\right)$ is not
of half-integer surgery type since it has no vectors of square 2.

\
  
Converting the half-integer surgery description above 
to integer surgery in the standard way gives a
cobordism $W$ from $\Sigma(J)$ to $\Sigma(K)$ whose intersection form
$Q_W$ is of half-integer surgery type.  Since $J$ is slice, $\Sigma(J)$ bounds a rational
homology ball $B$.  Joining $B$ to $W$ along $\Sigma(J)$ gives a smooth
closed four-manifold $X$ bounded
by $\Sigma(K)$.  The second Betti number of $X$ is twice the number of crossing changes used
to get from $K$ to $J$.

Suppose now that $K$ is converted to $J$ by changing $p$ positive and $n$ negative crossings,
with $n=\sigma(K)/2$.  Then $K$ bounds a disk in $B^4\#^{p+n}\cc\pp^2$.  Let
$X'$ be the double cover of the blown-up four-ball branched along this disk.
It follows from a theorem of Cochran and Lickorish \cite[Theorem 3.7]{cl} that 
$X'$ is positive-definite with $b_2(X')=2(p+n)$.  The following theorem is
based on the idea that in fact $X$ is diffeomorphic to $X'$.

\begin{maintheorem}
\label{thm:mainthm} Suppose that a knot $K$ may be converted to a slice knot by
changing $p$ positive and $n$ negative crossings, with
$n=\sigma(K)/2$. Then the branched double cover
$\Sigma(K)$ bounds a positive-definite smooth four-manifold $X$ with $b_2(X)=2(p+n)$
whose intersection form $Q_X$ is of half-integer surgery type, with exactly $n$ of 
$Q_X(x_1,x_1),\dots,Q_X(x_{p+n},x_{p+n})$ even,
and $\det Q_X$ divides $\det K$ with quotient a square.
\end{maintheorem}

For knots whose determinant is square-free, it follows that the first two 
parts of \ozsvath\ and \szabo's
obstruction to unknotting number one \cite[Theorem 1.1]{osu1} 
(without the symmetry condition) in fact give 
an obstruction to $u_s(K)=1$.

\begin{maincor}
\label{cor:os}
The knots
$$7_4, 8_{16}, 9_5, 9_{15}, 9_{17}, 9_{31}, 10_{19}, 
10_{20}, 10_{24}, 10_{36},  
10_{68}, 10_{69}, 10_{86},$$ 
$$10_{97}, 10_{105}, 
10_{109}, 10_{116}, 10_{121}, 10_{122}, 
10_{144}, 10_{163}, 10_{165}$$
have slice genus 1 and slicing number 2.
\end{maincor}

(Note that as mentioned above this was shown for $7_4$ in \cite{liv} and for $8_{16}$ in \cite{my}.
The slice genus information in Corollaries \ref{cor:os} and \ref{cor:u} is taken from \cite{knotinfo}.)

Furthermore, the obstruction given in \cite[Theorem 5]{u} to unknotting a knot by changing
$p$ positive and $n=\sigma(K)/2$ negative crossings is in fact an obstruction to slicing, provided
again that $\det K$ is square-free.  

\begin{maincor}
\label{cor:u} The knots 
$$9_{10}, 9_{13}, 9_{38}, 10_{53}, 10_{101}, 10_{120}$$ 
have slice genus $2$ and slicing number $3$.

The 11-crossing two-bridge knot $S(51,35)$ (Dowker-Thistlethwaite name $11a365$)
has slice genus $3$ and slicing number $4$.
\end{maincor}

It may also be shown that for some of the knots in Corollaries \ref{cor:os} and \ref{cor:u},
Livingston's invariant $U_s$ is not equal to the slice genus.  The knot $7_4$ is such an example (as Livingston suggested in \cite{liv}),
and in fact we find that it is the first member of an infinite family of such examples.

\begin{maincor}
\label{cor:Kn}
For each positive integer $n$,
there exists a two-bridge knot $K_n$ with signature $2n$ and slice genus $n$ which
cannot be sliced by changing $n$ negative crossings and any number of positive crossings;
hence $U_s(K_n)>n$.
\end{maincor}

\vskip2mm \noindent{\bf Acknowledgements.}
It is a pleasure to thank Andr\'{a}s Stipsicz and Tom Mark for helpful conversations.


\section{Proof of Theorem \ref{thm:mainthm}}
\label{sec:proof}

In this section we prove our main result.  

Recall that
a positive-definite integer-valued symmetric bilinear form $Q$ on a free abelian group $A$ gives an
integer lattice $L$ in Euclidean space on tensoring with $\rr$.  We say a lattice $L$ in
$\rr^n$ is of half-integer surgery type if the corresponding form $Q$ is (see Definition
\ref{def:halfint}).  Also a matrix representative for $Q$ is referred to as a Gram matrix for $L$.  For convenience we will frequently denote $Q(x,y)$ by $x\cdot y$, and $Q(x,x)$ by
$x^2$.

The proof of Theorem \ref{thm:mainthm} consists of a topological and
an algebraic step.  Following \cite{u} we show using careful analysis of Montesinos'
trick that, under the hypotheses of the theorem, $\Sigma(K)$ bounds
a positive-definite manifold $X$ and that the intersection
pairing of $X$ is of half-integer surgery type when restricted to some finite index sublattice.
We then show that if a lattice $M$ has an odd index sublattice $L$ of half-integer surgery
type then in fact $M$ is of half-integer surgery
type.

We begin with a couple of lemmas.

\begin{lemma}
\label{lem:mod4}
Let $Q$ be a form of half-integer surgery type, with $m_i=Q(y_i,y_i)$.  Then 
$$\det Q\equiv\prod_{i=1}^r(2m_i-1)\pmod4.$$
\end{lemma}
\proof This follows from the discussion after Lemma 2.2 in \cite{u}.\endproof

\begin{lemma}
\label{lem:mod2}
Let $Q$ be a block matrix of $r\times r$ blocks of the form
$\left(\begin{matrix}2I & *\\ * & *\end{matrix}\right)$
which is congruent modulo 2 to
$\left(\begin{matrix}2I & I\\I & X\end{matrix}\right)$.
Then there exists $P=\left(\begin{matrix}I & *\\0 & R\end{matrix}\right)\in GL(2r,\zz)$
with
$$P^T Q P=\left(\begin{matrix}2I & I\\I & X'\end{matrix}\right),$$
and $X'\equiv X\pmod2$.
\end{lemma}
\proof
Let $Q$ be the Gram matrix of a lattice with basis $x_1,\dots,x_r,z_1,\dots,z_r$. 
By successively adding multiples of $x_i$ to each of $z_1,\dots,z_r$ we get a new
basis $x_1,\dots,x_r$, $y_1,\dots,y_r$
with $x_i\cdot y_j=\delta_{ij}$; since $x_i$ has even square this preserves parities
on the diagonal.
\endproof

The following was originally proved by \ozsvath\ and \szabo\ \cite{osu1} in the case
$p+n=1$ and $J$ is the unknot.

\begin{proposition}
\label{prop:monttrick}
Suppose that a knot $K$ may be converted to a slice knot $J$ by
changing $p$ positive and $n$ negative crossings, with
$n=\sigma(K)/2$. Then the branched double cover
$\Sigma(K)$ bounds a positive-definite four-manifold $X$ with $b_2(X)=2(p+n)$.
The lattice $(H_2(X;\zz),Q_X)$ contains a finite index sublattice of half-integer type,
which has a basis as in Definition \ref{def:halfint} with exactly
$p$ elements of odd square.
\end{proposition}

\proof
We adapt the proof of \cite[Lemma 3.2]{u}.  By Montesinos' lemma 
(\cite{mont}, or see \cite[Lemma 3.1]{u}), $\Sigma(J)$ is the result
of surgery on some link $L$ in $S^3$ with half-integer framing coefficients.
Convert to integer surgery (see \cite{gs} or \cite[Lemma 2.2]{u}), and let
$Q_J$ be the resulting linking matrix of half-integer type.  We may assume
(after possibly adding a $-1/2$ framed unknot to $L$) that $\det Q_J$ is positive.
Denote by $X_J$ the two-handlebody with boundary $\Sigma_J$ and
intersection form represented by $Q_J$.

Note that since $J$ is slice it has signature zero and determinant $\det J=\det Q_J=k^2$ for some 
odd integer $k$.  Suppose that $K_-$ is a knot of signature 2 which may be converted
to $J$ by changing a single negative crossing $c$.  Then $\Sigma(K)$ is the result of surgery
on $L\cup C$ for some knot $C$ in $S^3$, with framing $(2m-1)/2$ on $C$.  Let $K_0$ be the
result of taking the oriented resolution of the crossing $c$; then as in \cite[Lemma 3.2]{u}
we have that $\Sigma(K_0)$ is surgery on $L\cup C$ with framing $m$ on $C$.  Converting to integer
surgery we find that $\Sigma(K_-),\Sigma(K_0)$ are given by integer surgeries with linking matrices
$$Q_-=\left(\begin{matrix}2&1&&0&\\1&m&&*&\\&&&&\\0&*&&Q_J&\\&&&&\end{matrix}\right),\quad
Q_0=\left(\begin{matrix}m&&*&\\&&&\\ *&&Q_J&\\&&&\end{matrix}\right).$$

Let $\Delta_K(t)$ denote the Conway-normalised Alexander polynomial of a knot $K$.  This satisfies
$$\Delta_K(-1)=(-1)^{\sigma(K)/2}\det K,$$
and thus
$$\Delta_J(-1)=k^2,\ \Delta_{K_-}(-1)=-|\det Q_-|=-|2\det Q_0 -k^2|.$$

The skein relation for the Alexander polynomial (see \cite{lick2}) then yields
$$k^2+|2\det Q_0 -k^2|=2|\det Q_0|,$$
from which we conclude that $\det Q_-=2\det Q_0 -k^2$ is positive.  Now using Lemma \ref{lem:mod4}
we have
$$2m-1\equiv\det Q_-=\det K_-\equiv3\pmod4,$$
and thus $m$ is even.
(The last congruence is due to Murasugi
\cite{m}: $\det K\equiv\sigma(K)+1\pmod4$ .)

Now suppose $K_+$ is a knot of signature 0 which may be converted to 
$J$ by changing a positive crossing.  Again $\Sigma(K)$ is
obtained by half-integer  surgery on $L\cup C$ for some knot $C$ in $S^3$, 
with framing $(2m-1)/2$ on $C$.  Let $Q_+$ denote the linking matrix after converting to
integer surgery.  A similar argument as above (see \cite{u} for the case where $J$ is the unknot)
shows that $\det Q_+>0$ and $m$ is odd.

Let $c_1,\dots,c_{p+n}$ be the set of crossings ($p$ positive, $n$ negative) 
in some chosen diagram of $K$ that we change to convert to $J$.
Then $\Sigma(K)$ is Dehn surgery on the link
$L\cup C_1\cup\dots\cup C_{p+n}$, with half-integer framing coefficients.  Each $C_i$
corresponds to a crossing $c_i$.
Dehn surgery on $L$ union a sublink of $C_1\cup\dots\cup C_{p+n}$ gives the double
branched cover of a knot which is obtained from $K$ by changing
a subset of the crossings $c_1,\dots,c_{p+n}$.  In particular surgery on
the knot $L\cup C_i$ yields the double branched cover of the knot
$K'_i$ which is obtained from $K$ by changing all of the crossings except $c_i$.
By the condition $n=\sigma(K)/2$ the knot signature changes every time
a negative crossing is changed and remains constant when a positive crossing changes.
It follows from the discussion above applied to
$K'_i$ that the framing on each $C_i$ is of the form $(2m_i-1)/2$ and exactly
those $m_i$ which
correspond to changing negative crossings of $K$ are even.

Denote by $X_K$ the two handlebody with boundary $\Sigma(K)$ that results 
from converting to integer surgery 
(i.e. surgery on the link $L\cup C_1\cup\dots\cup C_{p+n}$,
followed by surgery on a 2-framed meridian of each component).  Then
$X_K$ has intersection form of half-integer surgery type; moreover we can view
$X_K$ as the union of $X_J$ and a surgery cobordism $W$ along the common boundary
$\Sigma(J)$.  We will show by an induction argument 
that this cobordism is positive-definite.

Let $K_j$ denote the knot obtained from $K$ by changing crossings 
$c_{j+1},\dots,c_{p+n}$, and
let $Q_j$ be the linking matrix of half-integer type obtained 
by converting the corresponding
Dehn surgery diagram of $\Sigma(K_j)$ to integer type.  

Suppose that $\det Q_{j-1}$ is positive and hence equal to $\det K_{j-1}$.  We have
$$Q_j=\left(\begin{matrix}2&1&&0&\\1&m_j&&*&\\&&&&\\0&*&&Q_{j-1}&\\&&&&\end{matrix}\right),$$
and by Lemma \ref{lem:mod4}
$$\det Q_j\equiv(2m_j-1)\det Q_{j-1}\pmod4.$$
If $c_j$ is a positive crossing then $m_j$ is odd and so
\begin{equation}
\label{eqn:detQ}\det Q_j\equiv\det Q_{j-1}\pmod4.
\end{equation}
On the other hand
the signature of $K_{j-1}$ is equal to that of $K_j$ and so 
\begin{equation}
\label{eqn:detK}\det K_j\equiv\det K_{j-1}\pmod4.
\end{equation}
Comparing (\ref{eqn:detQ}) and (\ref{eqn:detK}) we see that $\det Q_j$ is congruent modulo
4 to its absolute value.  Since it is an odd number it must be positive.

On the other hand if $c_j$ is negative we find both congruences (\ref{eqn:detQ}) and (\ref{eqn:detK})
do \emph{not} hold, and again it follows that $\det Q_j$ is positive.

By induction we see that $Q_J=Q_0, Q_1, \dots ,Q_{p+n}$ all have positive determinants.  
Thus the surgery cobordism $W$ is built by attaching $2(p+n)$ two-handles to
the two-handlebody $X_J$, and before and after each attachment we have
a four-manifold whose intersection pairing has positive determinant.
It follows that $W$ is positive-definite.

We claim $(H_2(W;\zz),Q_W)$ contains a finite index sublattice with a 
basis as in Definition \ref{def:halfint} with exactly $p$ elements of odd square.
Suppose that $L$ has $r$ components, so that $b_2(X_K)=2(p+n+r)$.
Let $\{x_i,y_i\}$ be a basis for $H_2(X_K;\zz)$ as in Definition \ref{def:halfint},
chosen so that $\{x_i,y_i\}_{i>p+n}$ is a basis for the sublattice $H_2(X_J;\zz)$.
For $i\le p+n$, take $y_i$ to be the class corresponding to the two-handle
attached along $C_i$, and $x_i$ that corresponding to the two-handle
attached along the meridian of $C_i$.  Note that each $x_i$ is contained in
$H_2(W;\zz)$.  The rough idea is to form a sublattice by projecting the span
of $\{x_i,y_i\}_{i\le p+n}$ orthogonally to $H_2(W;\zz)$.  Of course we
cannot quite do this; however since $H_1(\Sigma(J);\zz)$ has order $k^2$
we may write
$$k^2 y_i=z_i+w_i,\qquad i=1,\dots,p+n$$
with $z_i\in H_2(W;\zz)$ and $w_i\in H_2(X_J;\zz)$.
We claim that the self-intersection of $z_i$ has the same parity as that of $y_i$.
To see this note that for each $j>p+n$, $y_i$ is orthogonal to $x_j$ and
hence so is $w_i$.  It follows that $w_i$ is in the span of $\{x_i\}_{i>p+n}$
and so has even self-intersection.
Thus we have a full rank sublattice of $(H_2(W;\zz),Q_W)$ with basis
$\{x_i,z_i\}_{i\le p+n}$; by Lemma \ref{lem:mod2}, this sublattice
has a basis as in Definition \ref{def:halfint} with exactly $p$ 
elements of odd square.

Form the manifold $X$ with one boundary component by capping off the $\Sigma(J)$ end of
$W$ with the rational ball $B$ given as the double branched cover of $B^4$ along a slice disk
bounded by $J$ \cite{cg}.
Then $(H_2(W;\zz),Q_W)$ is a sublattice of $(H_2(X;\zz),Q_X)$, and therefore so is 
$\Span(\{x_i,z_i\}_{i\le p+n})$.
\endproof

\begin{proposition}
\label{prop:subhalfint}
Suppose $M$ is a positive-definite integer lattice of rank $2r$ and $L$ is a sublattice
of $M$ of odd index $l$.  If $L$ is of half-integer type then so is $M$.  Moreover,
the number of elements of a basis of $M$ as in Definition \ref{def:halfint} with
odd square is the same as that for $L$.
\end{proposition}

This will follow from the next two lemmas, the first of which is standard.

\begin{lemma}
\label{lem:gln}
The natural map from $GL(n,\zz)$ to $GL(n,\zz/2)$ is onto for any $n$.
\end{lemma}

\proof Use induction on $n$.  Suppose that $R\in M(n,\zz)$ has odd determinant.
The cofactor expansion across the first row yields
$$\det R=r_{11}R_{11}+r_{12}R_{12}+\cdots+r_{1n}R_{1n}.$$
Since the determinant is odd, so is at least one $r_{1j}R_{1j}$.
By induction we may choose $\tilde{R}\equiv R\pmod2$ with
$\tilde{R}_{1j}=1$, then adjust the value of $r_{1j}$ to get
$\det\tilde{R}=1$.
\endproof

\begin{lemma}
\label{lem:mod2basis}
Suppose $M$ is a positive-definite integer lattice of rank $2r$, and that $L$ 
is a lattice of half-integer type which is a sublattice
of $M$ of odd index $l$.
Let $x_1,\dots,x_r,y_1,\dots,y_r$ be
a basis for $L$ as in Definition \ref{def:halfint}, and let $Q_L$ be the Gram
matrix of $L$ in this basis.  Then $x_1,\dots,x_r$ may be extended to a basis
$x_1,\dots,x_r,z_1,\dots,z_r$ for $M$ with
$$Q_M\equiv Q_L\pmod2.$$
\end{lemma}

\proof
Let $m_i=y_i\cdot y_i$.
In the given basis $Q_L$ is in block form
$\left(\begin{matrix}2I & I\\I & X\end{matrix}\right)$, with $\diag(X)=(m_1,\dots,m_r)$.
By Theorem 6 in Chapter 1 \S 3 of \cite{grublek} there exists a basis
$a_1,\dots,a_r,z_1,\dots,z_r$ of $M$ with $x_i\in\Span_\zz\{a_1,\dots,a_i\}$.
A simple induction argument using the fact that $x_i\cdot y_j=\delta_{ij}$ shows that
in fact (after possibly multiplying by $-1$) we have $a_i=x_i$ for $i=1,\dots,r$.

Let $P\in M(2r,\zz)$ be the matrix whose $i$th column is the coefficient vector of the
$i$th basis vector of $L$ in the basis $x_1,\dots,x_r,z_1,\dots,z_r$.  Then $P$
is in block form 
$\left(\begin{matrix}I & *\\0 & R\end{matrix}\right)$,
and 
$$Q_L=P^T Q_M P.$$
Since $\det Q_L=l^2\det Q_M$ we have $\det R=l$ is odd.  By Lemma \ref{lem:gln} we may choose
$\tilde{R}\in GL(r,\zz)$ with $\tilde{R}\equiv R\pmod2$.  Applying the transition
matrix $\tilde{P}=\left(\begin{matrix}I & *\\0 & \tilde{R}\end{matrix}\right)$
to $M$ yields the required basis.
\endproof

\noindent{\it Proof of Proposition \ref{prop:subhalfint}.}
Let $x_1,\dots,x_r,z_1,\dots,z_r$ be the basis of $M$ given by Lemma \ref{lem:mod2basis},
in which
$$Q_M=\left(\begin{matrix}2I & *\\ * & *\end{matrix}\right)
\equiv Q_L=\left(\begin{matrix}2I & I\\I & X\end{matrix}\right)\pmod2.$$
The proposition now follows from Lemma \ref{lem:mod2}.
\endproof

\noindent{\it Proof of Theorem \ref{thm:mainthm}.}
Let $X$ be the four-manifold bounded by $\Sigma(K)$ given by Proposition \ref{prop:monttrick}.
The order of the first homology of a knot is always odd; it follows that $Q_X$ has odd order and
the sublattice given by Proposition \ref{prop:monttrick} has odd index in $(H_2(X;\zz),Q_X)$.
 Theorem \ref{thm:mainthm} now follows immediately from Proposition
\ref{prop:subhalfint}.\endproof


\section{Examples}
\label{sec:examples}
Theorem \ref{thm:mainthm}
tells us that to show that a knot $K$ cannot be converted to a slice knot by changing
$p$ positive and $n=\sigma(K)/2$ negative crossings, we
must show that its double branched cover 
cannot bound a four-manifold with an intersection form with certain 
properties.

We will make use of two very effective gauge-theoretic 
obstructions to a rational 
homology three-sphere $Y$ bounding
a positive-definite form $Q$.  On the one hand, \ozsvath\ and \szabo\ define
a function
$$d:\spinc(Y)\to\qq$$
coming from the absolute grading in Heegaard Floer homology, and they show that for
each \spinc structure $\spincs$ on a positive-definite four manifold bounded by $Y$
the following must hold:
\begin{eqnarray}
c_1(\spincs)^2-b_2(X)&\ge& 4d(\spincs|_Y),\label{eqn:ineq}\\
\mbox{and}\quad
c_1(\spincs)^2-b_2(X)&\equiv& 4d(\spincs|_Y) \pmod2.\label{eqn:cong}
\end{eqnarray}
The left hand side depends on the intersection form of $X$.

To determine if a given knot $K$ may be sliced by changing $p$ positive and 
$n=\sigma(K)/2$ negative crossings we carry out the following steps:
\renewcommand{\labelenumi}{\arabic{enumi}.}
\begin{enumerate}
\item Compute $d:\Sigma(K)\to\qq$;
\item Find a complete finite set of representatives $Q_1,\dots,Q_m$ of forms
of rank $2(p+n)$ satisfying the conclusion of Theorem \ref{thm:mainthm};
\item Check using (\ref{eqn:ineq},\ref{eqn:cong}) whether $\Sigma(K)$ is obstructed from
bounding each of $Q_1,\dots,Q_m$.
\end{enumerate}
Details on Heegaard
Floer theory and the $d$ invariant may be found in \cite{os4,os6,osu1}; 
for a summary of
how this theory may be used in our context see \cite{u}.

Another approach to understanding the set of positive-definite forms
that a three-manifold $Y$ may bound is to make use of 
Donaldson's diagonalisation theorem \cite{d}.
This approach works well for Seifert fibred rational homology spheres, and especially
for lens spaces (see e.g.~ \cite{lisca}, \cite{det3}).
Knowledge of a particular negative-definite four-manifold $X_1$ bounded
by $Y$ constrains the intersection form of a positive-definite $X_2$ with the
same boundary, since $X=X_2\cup-X_1$ is a closed positive-definite manifold
with $(H_2(X;\zz),Q_X)\cong\zz^m$ for some $m$.  We will illustrate this technique
in Subsection \ref{subsec:Kn}.

The slicing number of a knot is the same as that for its 
reflection.  We assume in what follows
that all knots have nonnegative signature.  (This distinguishes between the knot and
its reflection unless the signature is zero.)

\subsection{Knots with slice genus one}
For a knot with $\sigma(K)=2$ and $u_s(K)=1$, it follows from 
inequality (\ref{eqn:nsig}) and 
Theorem \ref{thm:mainthm} that $\Sigma(K)$ bounds a four-manifold whose intersection
form is represented by the matrix
$$\left(\begin{matrix}m&1\\1&2\end{matrix}\right),$$
with $\det K=(2m-1)t^2$ for some integer $t$.

For a knot with $\sigma(K)=0$ and $u_s(K)=1$ we find that either $\Sigma(K)$ or
$-\Sigma(K)$ must bound such a positive-definite four-manifold.

The knots listed in Corollary \ref{cor:os} have square-free determinant.  For each
of them, \ozsvath\ and \szabo\ have shown in \cite{osu1}, 
using (\ref{eqn:ineq},\ref{eqn:cong}), that $\pm\Sigma(K)$ cannot
bound
$$\left(\begin{matrix}\frac{\det K +1}{2}&1\\1&2\end{matrix}\right);$$
it is also known in each case that the knot can be
unknotted with two crossing changes.  We conclude that $u_s(K)=2$.

\begin{remark}
Each of the knots $10_{29}$, $10_{40}$, $10_{65}$, $10_{67}$, 
$10_{89}$, $10_{106}$, and $10_{108}$ has signature 2 and 
$\det K=st^2$ for some $t>1$.  In each case
\ozsvath\ and \szabo\ have shown $\Sigma(K)$ cannot bound
$\left(\begin{matrix}\frac{\det K +1}{2}&1\\1&2\end{matrix}\right)$
and hence $K$ has unknotting number two.

However we find in each case $\Sigma(K)$  is \emph{not} obstructed from bounding
$\ds\left(\begin{matrix}\frac{s+1}{2}&1\\1&2\end{matrix}\right)$.
\end{remark}


\subsection{Knots with slice genus two or three}
We now consider the knots in Corollary \ref{cor:u}.  Each of 
$9_{10}, 9_{13}, 9_{38}, 10_{53}, 10_{101}, 10_{120}$
has signature 4 and slice genus 2.  In \cite{u} we have shown, 
using (\ref{eqn:ineq},\ref{eqn:cong}), that
in each case $\Sigma(K)$ cannot bound a positive-definite form $Q$ as in Definition 
\ref{def:halfint}
with rank 4, $\det Q=\det K$ and $Q(x_i,x_i)$ even.
Since the knots have square-free determinant, it follows from Theorem \ref{thm:mainthm}
that they cannot be sliced with two crossing changes.  Each can be unknotted with 
three crossing changes, and so in each case $u_s(K)=3$.

Similarly $K=11a365$ is shown in \cite{u} to have unknotting number 4, and since its
determinant is square-free the same argument shows it has $u_s(K)=4$.

\subsection{Knots with large slice genus}
\label{subsec:Kn}
\begin{figure}[htbp]
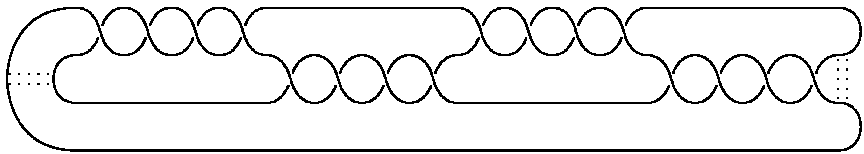

\begin{center}
\ifpic
\leavevmode
\xygraph{
!{0;/r1.5pc/:}
!{\hcap[-3]}[d]
!{\hcap[-1]}[u]
!{\htwist}[ldd]
!{\xcaph[1]@(0)}[ld]
!{\xcaph[1]@(0)}[uuu]
!{\htwist}[ldd]
!{\xcaph[1]@(0)}[ld]
!{\xcaph[1]@(0)}[uuu]
!{\htwist}[ldd]
!{\xcaph[1]@(0)}[ld]
!{\xcaph[1]@(0)}[uuu]
!{\htwist}[ldd]
!{\xcaph[1]@(0)}[ld]
!{\xcaph[1]@(0)}[uuu]
!{\xcaph[1]@(0)}[ld]
!{\htwist}[ldd]
!{\xcaph[1]@(0)}[uuu]
!{\xcaph[1]@(0)}[ld]
!{\htwist}[ldd]
!{\xcaph[1]@(0)}[uuu]
!{\xcaph[1]@(0)}[ld]
!{\htwist}[ldd]
!{\xcaph[1]@(0)}[uuu]
!{\xcaph[1]@(0)}[ld]
!{\htwist}[ldd]
!{\xcaph[1]@(0)}[uuu]
!{\htwist}[ldd]
!{\xcaph[1]@(0)}[ld]
!{\xcaph[1]@(0)}[uuu]
!{\htwist}[ldd]
!{\xcaph[1]@(0)}[ld]
!{\xcaph[1]@(0)}[uuu]
!{\htwist}[ldd]
!{\xcaph[1]@(0)}[ld]
!{\xcaph[1]@(0)}[uuu]
!{\htwist}[ldd]
!{\xcaph[1]@(0)}[ld]
!{\xcaph[1]@(0)}[uuu]
!{\xcaph[1]@(0)}[ld]
!{\htwist}[ldd]
!{\xcaph[1]@(0)}[uuu]
!{\xcaph[1]@(0)}[ld]
!{\htwist}[ldd]
!{\xcaph[1]@(0)}[uuu]
!{\xcaph[1]@(0)}[ld]
!{\htwist}[ldd]
!{\xcaph[1]@(0)}[uuu]
!{\xcaph[1]@(0)}[ld]
!{\htwist}[ldd]
!{\xcaph[1]@(0)}[uuu]
!{\hcap}[dd]
!{\hcap}[u(0.9)]
!{\knotstyle{.}}
!{\xcapv[0.8]@(0)}[u(1)r(0.2)] 
!{\xcapv[0.8]@(0)}[u(0.7)l(17.6)] 
!{\xcaph[0.8]@(0)}[l(1)d(0.2)]
!{\xcaph[0.8]@(0)}
}
\else \vskip 5cm \fi
\begin{narrow}{0.3in}{0.3in}
\caption{
{\bf{The knot $K_2$.}}
The two pairs of dashed arcs indicate where to attach ribbons to go from $K_n$ to $K_{n-1}$.}
\label{fig:Kn}
\end{narrow}
\end{center}
\end{figure}

In this subsection we prove Corollary \ref{cor:Kn}.  We define $K_n$ to be the 4-plat closure of
the four-string braid $(\sigma_1{}\!^{4}\sigma_2{}\!^{4})^n$, 
as illustrated in Figure \ref{fig:Kn} for
$n=2$.  For $n=1$ this is the knot 
$7_4=S(15,11)$ shown by Lickorish to have unknotting number 2
and by Livingston to have slicing number 2 \cite{lick,liv}.

As illustrated in the diagram, two oriented ribbon moves convert $K_n$ to $K_{n-1}$.  Since $K_0$
is the unknot this shows that the slice genus of $K_n$ is at most $n$.  The signature of
$K_n$ may be shown (see below) to be $2n$.  
We conclude that 
$$g_s(K_n)=n.$$

Let $P(a_1,a_2,\dots,a_m)$ denote the plumbing of
disk bundles over two-spheres corresponding to the linear graph with $m$ vertices, where
the $i$th vertex has weight $a_i$.
The double cover of $S^3$ branched along $K_n$ is the boundary of $P(4,4,\dots,4)$.
Let $Q_n$ denote the rank $2n$ intersection form of this plumbing,
and let $L_n$ denote the associated lattice in $\rr^{2n}$.
The following
sequence of blow-ups and blow-downs exhibits $\Sigma(K_n)$ as the boundary of a negative-definite
plumbing:

\begin{eqnarray*}
\Sigma(K_n)&\cong&\partial P(4,4,\dots,4)\\
&\cong&\partial P(-1,2,-1,2,\dots,-1,2,-1)\\
&\cong&\partial P(-2,-1,1,-2,-1,1,\dots,-2,-1,1,-1)\\
&\cong&\partial P(-2,-2,-3,-2,-3,-2,\dots,-3,-2,-3,-2,-2).\\
\end{eqnarray*}

(We note the above shows that $K_n$ may also be represented by the
alternating diagram which is the closure of the braid
$$(\sigma_1{}\!^{-1}\sigma_2{}\!^{2})^{2n}\sigma_1{}\!^{-1};$$
using the formula of Gordon and Litherland \cite{gl}
it is easy to compute the signature from this diagram.)

Let $X'_n$ denote the positive-definite plumbing $P(2,2,3,2,3,2,\dots,3,2,2)$
whose boundary is $-\Sigma(K_n)$.  Let $Q'_n$ denote its intersection form
and let $L'_n$ denote the associated lattice.
Note that $L'_n$ has dimension $2n+3$: there are $n$ vertices
with weight $3$ and $n+3$ with weight $2$.

\begin{lemma}
\label{lem:Knobst1}
Suppose $\Sigma(K_n)$ is given as the boundary of a smooth four-manifold
$X$ with positive-definite intersection form $Q_X$.  Then
$(H_2(X;\zz),Q_X)$ embeds as a finite index sublattice of
$L_n\oplus \zz^k$
for some $k\ge0$.
\end{lemma}
\proof  Gluing $X$ to $X'_n$ along their boundary gives a 
closed positive-definite
manifold.   It follows from Donaldson's theorem
that the lattice $L'_n$ embeds as a sublattice of $\zz^m$
with the lattice $(H_2(X;\zz),Q_x)$ contained in its orthogonal complement.

Let $e_1,\dots,e_m$ be an orthonormal basis of $\zz^m$.
Up to automorphism of $\zz^m$ there is a unique way to embed
$L'_n$: the first vertex vector must map to $e_1+e_2$, the second to
$e_2+e_3$, the third to $e_3+e_4+e_5$, the fourth to $e_5+e_6$ and so on.
Thus the image of $L'_n$ is contained
in a $\zz^{3n+4}$ sublattice of $\zz^m$.
An easy calculation shows that the orthogonal complement of
$L'_n$ in $\zz^{3n+4}$ is spanned by the vectors
$e_1-e_2+e_3-e_4$, $e_4-e_5+e_6-e_7$,\dots, 
$e_{3n+1}-e_{3n+2}+e_{3n+3}-e_{3n+4}$.
These span a copy of $L_n$, from which the conclusion follows.\endproof

\begin{lemma}
\label{lem:Knobst2}
For any $n\ge1$, $k\ge0$,  the lattice $L_n\oplus \zz^k$
does not admit any finite index sublattices of half-integer surgery type.
\end{lemma}
\proof
For $k=0$ this is immediate since $L_n$ has no nonzero 
vectors of square less than $4$.
If $k>0$ let $e_1$,\dots,$e_k$ be an orthonormal basis of $\zz^k$, and suppose
we have a sublattice of $L_n\oplus \zz^k$ of half-integer surgery type
with basis $\{x_i,y_i\}$ as in Definition \ref{def:halfint}.
Up to an automorphism of $L_n\oplus \zz^k$ we have $x_1=e_1+e_2$.  Then
$x_2$ is orthogonal to $x_1$.  We cannot have $x_2=e_1-e_2$ since $y_1$ pairs
evenly with $x_2$ and oddly with $x_1$.  Thus up to automorphism, $x_2=e_3+e_4$.
It follows that any sublattice of $L_n\oplus \zz^k$ of half-integer surgery type
has rank at most $k$.\endproof

Corollary \ref{cor:Kn} now follows from Lemmas \ref{lem:Knobst1} and \ref{lem:Knobst2}
and Theorem \ref{thm:mainthm}.

\begin{remark}
\label{rem:Kn}
Livingston conjectured in \cite{liv} that the difference 
$U_s-g_s$ can be arbitrarily large.
It is possible to unknot $K_n$ by changing $2n$ positive crossings
in the diagram as in Figure \ref{fig:Kn}.
In the absence of any further evidence it is 
tempting to conjecture that
for these knots $U_s-g_s=n$; in any case this would seem to be a good candidate 
with which to attempt to verify Livingston's conjecture.
\end{remark}

\begin{remark}
 \label{rem:disk}
The trace of a homotopy from a knot $K$ to a slice knot $J$ is an immersed annulus in $S^3\times I$; capping this off with a slice disk yields an immersed disk $D$ in $B^4$ bounded by $K$.  If $J$ is obtained from $K$ by changing $p$ positive and $n$ negative crossings then the resulting disk $D$ has $p$ negative self-intersections and
$n$ positive self-intersections.  (This is why changing a positive crossing is often referred to as a ``negative crossing change''.)
Instead of considering crossing changes one may ask whether a knot $K$ bounds an immersed disk in $B^4$ with a prescribed number of self-intersections, or with prescribed numbers of self-intersections of each sign.

Rudolph has shown in \cite{r} that the minimal number of self-intersections in a \emph{ribbon} immersed disk bounded by $K$ is equal to the minimal number of crossing changes to get from $K$ to a ribbon knot.  (Here a ribbon surface in $B^4$ is one on which the radial distance function has no maxima, and a ribbon knot is a knot which bounds an embedded ribbon disk.)  Knowing whether a result analagous to Rudolph's for non-ribbon disks and slice knots holds would be very interesting.  It is to be expected that the conclusion of Theorem \ref{thm:mainthm} holds under the weaker hypothesis that $K$ bounds an immersed disk with $p$ negative, $n=\sigma(K)/2$ positive self-intersections.  The expected proof would generalise that of \cite{my} and also make use of \cite[Theorem 3.7]{cl}.

\end{remark}


\end{document}